\documentclass[a4paper, 11pt]{article}	 

\usepackage[paper=a4paper,left=1.8cm,right=1.8cm]{geometry}
\usepackage{cite}
\usepackage{colortbl}
\usepackage{color}
\usepackage{footnote}
\usepackage{subfig}
\usepackage{rotating} 
\usepackage[onehalfspacing]{setspace}
\usepackage{amsmath}
\usepackage{amsthm}
\usepackage{multirow}
\usepackage{amssymb}
\usepackage{psfrag}
\usepackage{graphicx}
\usepackage{algorithm}
\usepackage{algpseudocode}
\usepackage[T1]{fontenc}
\usepackage[utf8]{inputenc}
\usepackage{authblk}

\title{\LARGE \bf Simulation studies on online constraint removal with a Lyapunov function}
\author[1]{Michael Jost}
\author[2]{Gabriele Pannocchia}
\author[1]{Martin M\"onnigmann\thanks{Corresponding author. Email: martin.moennigmann@rub.de}}
\affil[1]{{\small{\em Automatic Control and Systems Theory, Ruhr-Universit\"at Bochum, Bochum, Germany}}}
\affil[2]{{\small{\em Department of Civil and Industrial Engineering, University of Pisa, Pisa, Italy}}}

\newcommand{\R}{{\mathbb R}}

\newcommand{\allIndices}{{\mathcal Q}}
\newcommand{\activeSet}{{\mathcal A}}
\newcommand{\inactiveSet}{{\mathcal I}}
\newcommand{\DOA}{{\mathcal{X}}}

\newcommand{\condNum}[1]{\kappa(#1)}

\begin{document}
\maketitle
\thispagestyle{plain}
\pagestyle{plain}
\section{Introduction}
We apply the method proposed in \cite{Jost2014} to 36 MPC implementations, which result from combining 
six sample receding horizon control problems (see Tab.~\ref{tab:systems}) with six QP solvers (see Tab.~\ref{tab:SummaryOptAlgs}). We implement each of the 36 system-solver-combinations both with and without constraint removal and compare computational times for statistically relevant numbers of runs. 

\section{Brief problem statement}\label{sec:mpc}
We introduce the notation only as needed to make this report self-contained and refer to \cite{Jost2014} and references therein for details. 
We consider discrete-time state space systems
\begin{align} \label{eq:stateSpaceSystem}
  x(t+1) &= Ax(t) + Bu(t),
\end{align}
with state $x(t)\in\R^n$, input $u(t)\in\R^m$ and matrices $A\in\R^{n\times n}$, $B\in\R^{n\times m}$, where $(A,B)$ is assumed to be stabilizable. The systems are subject to input and state constraints of the form
\begin{align}\label{eq:Constraints}
\begin{split}
u(t)\in\mathbb{U}\subset\R^m,\quad x(t)\in\mathbb{X}\subset\R^n,
  \end{split}
\end{align}
where $\mathbb{U}$ and $\mathbb{X}$ are polytopes (i.e.\ intersections of a finite number of halfspaces) that contain the origin in their interiors. 

The MPC problem for \eqref{eq:stateSpaceSystem}, \eqref{eq:Constraints} reads
\begin{align}\label{eq:MPCProblem}
\begin{array}{rrl}
&\min\limits_{U,X} &  x(N)^\prime P x(N) + \sum\limits_{k=0}^{N-1}  \left( x^\prime(k) Q x(k) +  u^\prime(k) R u(k)\right)\\ 
\mbox{s.\ t.}&
  x(k+1) &= A x(k) + B u(k), \; k = 0, \ldots, N-1\\
  &x(0) &= x_0 , \\
  &x(k) & \in \mathbb{X} , \; k= 1, \dots, N-1,\\
  &x(N) &\in \mathbb{X}_f , \\
 &u(k) &\in \mathbb{U} , \; k= 0, \dots, N-1,
\end{array}
\end{align}
where $U= \left(u^\prime(0), \dots, u^\prime(N-1)\right)^\prime$ and $X= \left(x^\prime(1), \dots, x^\prime(N)\right)^\prime$, $P\in\R^{n\times n}$, $P\succeq 0$,  $Q\in\R^{n\times n}$, $Q\succeq0$ and $R\in\R^{m\times m}$, $R\succ0$ and where $x_0$ is the initial condition. We assume $\mathbb{X}_f \subseteq \mathbb{X}$ to be a polyhedral terminal set that contains the origin in its interior. By a slight abuse of notation the initial condition is denoted by $x$ instead of $x_0$ in the remainder of the paper.

The MPC problem \eqref{eq:MPCProblem} can equivalently be stated in the form 
\begin{align}\label{eq:ParametricMPC}
\begin{split}
   \min\limits_{U} & \quad V(x,U) , \\ 
   \mbox{s.\ t.}\;\; & GU -w -E x \le 0
\end{split}
\end{align}
with cost function $V(x,U) = \frac{1}{2}x^\prime Y x + U^\prime F x  + \frac{1}{2}U^\prime H U$ and 
$H\in\R^{mN\times mN}$, $Y\in\R^{n\times n}$, $F\in\R^{mN\times n}$, $G\in\R^{q\times mN}$, $E\in\R^{q\times n}$ and $w\in\R^{q}$. The number of inequality constraints is denoted by $q$, and $\allIndices=\left\{1,\cdots,q\right\}$ is the index set of all constraint indices. For later use we note that the solution of the corresponding {\it unconstrained} optimization problem, i.e.~\eqref{eq:ParametricMPC} where all constraints are dropped, is given by
\begin{align}\label{eq:unconstrainedSolution}
 U^\star(x) = -H^{-1}F^\prime x.
\end{align}

\subsection{Notation}
Let $G^i$ denote the $i$-th row of $G$. For any ordered index set $\mathcal{W}\subset\allIndices$, let and $G^{\mathcal{W}}$ refer to the submatrix of $G$ that contains the rows indicated by $\mathcal{W}$. For all other matrices, the corresponding notation applies.  

Let $\DOA$ be the set of states for which the quadratic program~\eqref{eq:ParametricMPC} is feasible. We refer to the optimal solution of the quadratic program~\eqref{eq:ParametricMPC} by $U^\star$. More precisely, the function $U^\star:\DOA\rightarrow \mathbb{U}^N$ is defined by
\begin{align*}
U^\star(x)\;\; =\;\;& \underset{U} {\mathrm{arg\, min}} \quad V(x,U)\\
   \mbox{s.\ t.}\;\; & GU -w -E x \le 0.
\end{align*}
Constraint $i$ is called active for the state $x\in\DOA$ (or just active, for short) if it holds with equality at the optimal solution $U^\star(x)$, i.e.\ $G^iU^\star(x)=w^i+E^ix$. It is called inactive otherwise. We define the index sets of active and inactive constraints 
\begin{align}\label{eq:DefInactiveAndActiveSet}
 \begin{split}
  \activeSet(x) = \left\{i\in \allIndices \left|G^iU^\star(x)=w^i+E^ix\right.\right\},\\
  \inactiveSet(x)= \left\{i\in \allIndices \left|G^iU^\star(x)<w^i+E^ix\right.\right\},
 \end{split}
\end{align}
respectively. 
Note that $\allIndices = \activeSet(x)\cup\inactiveSet(x)$ and $\activeSet(x)\cap\inactiveSet(x)=\emptyset$ follows from this definition.

\section{Constraint removal}\label{sec:method}
Essentially, we show in~\cite{Jost2014} how to find a subset of the inactive constraints $\mathcal{J}(x)\subseteq\inactiveSet(x)$ {\it before} actually solving the quadratic program \eqref{eq:ParametricMPC}. Inactive constraints can be removed from \eqref{eq:ParametricMPC}. The quadratic program~\eqref{eq:ParametricMPC} therefore can be simplified to 
\begin{align}\label{eq:ReducedCondensedQP}
  \begin{split}
      \min\limits_{U} & \quad\frac{1}{2}x^\prime Y x + U^\prime F x  + \frac{1}{2}U^\prime H U, \\ 
       \mbox{s.\ t.}\;\; &G^iU -w^i -E^i x \le 0, \quad i\in\allIndices\backslash\mathcal{J}(x).
  \end{split}
\end{align} 
The construction of the subset $\mathcal{J}(x)$ is based on the decent property of a Lyapunov function of the closed-loop system (the optimal cost function of \eqref{eq:ParametricMPC}) \cite{Jost2014}.  Since only simple arithmetic operations are required to determine $\mathcal{J}(x)$, the computational cost of MPC can be reduced by first determining $\mathcal{J}(x)$, and then solving the simplified, smaller QP \eqref{eq:ReducedCondensedQP} instead of \eqref{eq:ParametricMPC}. 
Algorithm~\ref{alg:constraintRemovalAlgorithm} summarizes the steps that need to be carried out online. 

\begin{algorithm}[b]
\begin{algorithmic}[1]
\State{\bf Input} $x$
\State Calculate the set of inactive constraints $\mathcal{J}(x)$. 
\If{${\mathcal{J}(x)}= \mathcal{Q}$}
\State QP is unconstrained: $U^\star(x) = -H^{-1}F^\prime x$
\Else
\State Set up reduced QP~\eqref{eq:ReducedCondensedQP}. 
\State Solve reduced QP~\eqref{eq:ReducedCondensedQP} for $U^\star(x)$.
\EndIf
\State {\bf Output:} $u^\star(x)=
[\begin{matrix}
 I^{m \times m}&0&\cdots&0
 \end{matrix} ]U^\star(x)$
\end{algorithmic} \caption{MPC with constraint removal.}\label{alg:constraintRemovalAlgorithm} 
\end{algorithm}

\section{Simulation study} \label{sec:simstudy}
Six MPC problems \eqref{eq:MPCProblem} serve as examples in this study (six \textit{examples}, for short). 
These six examples are constructed from four constrained linear discrete-time systems~\eqref{eq:stateSpaceSystem} (four \textit{systems}, for short), which are described in Sect.~\ref{subsec:examples} and summarized in Tab.~\ref{tab:systems}. Six examples result from four systems, because one of the systems (labeled MIMO) is considered for two different horizons (30 and 75), and both with and without an a-priori, offline removal of redundant constraints.
For each example we implement MPC with six different QP solvers for a total of 36 example-solver combinations (36 \textit{combinations}). The six QP solvers are described in Sect.~\ref{subsec:solvers} and summarized in Tab.~\ref{tab:SummaryOptAlgs}.
For each of the resulting 36 combinations we compare the case with constraint removal to that without constraint removal (72 \textit{cases}). The abbreviations \textit{CR-MPC} and \textit{full-MPC} refer to any implementation with and without constraint removal, respectively.

We generate random initial values $x\in\DOA$ for every system and calculate trajectories for the MPC-controlled system until $\|x(t)\|\le10^{-3}$ for every initial condition. The specific numbers of initial conditions and QPs are summarized in Tab.~\ref{tab:systems}. Note that at least several hundred thousand QPs are solved in every case. The same initial conditions are used in all cases for a given system. 

For each of the 36 combinations, we compare the computational times needed to find the control law with Alg.~\ref{alg:constraintRemovalAlgorithm}, i.e.~MPC with constraint removal, to the computational time required for solving the full quadratic program~\eqref{eq:ParametricMPC}, i.e.~MPC implemented with the same QP solver but without constraint removal. 

We stress that the computational times reported for the cases with constraint removal include the times 
needed to construct the set $\mathcal{J}(x)$ and to set up and solve the reduced quadratic program~\eqref{eq:ReducedCondensedQP}. 

\begin{table}[t]
\centering
\caption{Summary of the linear discrete-time state systems that serve as sample systems.}\label{tab:systems}
\begin{tabular}{c||c|c|c|c|c||c|c}
  sample system & $n$ & $m$ & $N$ & $mN$ & $q$ & \# $x_0$ & \# QPs \\ \hline \hline
   MIMO30 & 10 & 3 & 30 & 90 & 780 & 4935 & 810608 \\\hline
   MIMO75 & 10 & 3 & 75 & 225 & 1950&  5046 & 829790 \\\hline
   MIMORED30 & 10 & 3 & 30 & 90 & 556 & 4708 & 771678 \\\hline
   ACC25 & 4 & 1 & 25 & 25 & 258  & 6159 & 1211507 \\\hline
   INPE50 & 4 & 1 & 50 & 50 & 500 & 6764 & 1112827 \\\hline
   COMA40 &12 &3 &40 &120 & 1200 & 8028 & 1136855
 \end{tabular}
\end{table}

\subsection{Examples}\label{subsec:examples}
We state the most important features of the examples in this section.

\paragraph*{{\bf MIMO30:}} The system in this example is the zero-order hold discretization of the continuous-time transfer function
 \begin{align}\label{eq:example}
  G(s) = 
  \begin{pmatrix}
   \frac{-5s+1}{36s^2+6s+1} & \frac{0.5}{8s+1}           & 0\\
   0                        & \frac{0.1(-10s+1)}{s(8s+1)} & \frac{-0.1}{(64s^2+6s+1)s}\\
   \frac{-2s+1}{12s^2+3s+1} &  0& \frac{2(-5s+1)}{16s^2+2s+1}
  \end{pmatrix},
 \end{align}
 with sample time $T_s=1\text{s}$. After removing uncontrollable states from~\eqref{eq:example}, a state space model \eqref{eq:stateSpaceSystem} with $n= 10$ states, $m= 3$ inputs and system matrices

 {\tiny
\begin{align*}
A &= 
\begin{pmatrix}
  8.34e-01  &  1.81e-01  &  7.27e-02  &  -5.61e-02  &  -1.59e-02  &  4.28e-03  &  -1.95e-03  &  -6.74e-03  &  -5.56e-03  &  -7.94e-03 \\
 -1.02e-01  &  9.38e-01  &  -5.11e-03  &  -1.62e-01  &  -1.44e-02  &  2.39e-03  &  1.19e-03  &  -4.30e-03  &  2.66e-03  &  3.95e-03 \\
 -4.09e-02  &  1.37e-01  &  8.91e-01  &  3.14e-01  &  1.02e-02  &  9.60e-04  &  4.77e-04  &  -1.73e-03  &  1.07e-03  &  1.59e-03 \\
 3.16e-02  &  1.50e-02  &  -1.37e-01  &  8.62e-01  &  -1.55e-02  &  -7.41e-04  &  -3.68e-04  &  1.34e-03  &  -8.27e-04  &  -1.23e-03 \\
 8.94e-03  &  7.21e-03  &  -8.31e-03  &  4.17e-02  &  8.84e-01  &  -2.10e-04  &  -1.04e-04  &  3.78e-04  &  -2.34e-04  &  -3.47e-04 \\
 -2.41e-03  &  -2.00e-03  &  1.39e-03  &  1.30e-03  &  -1.74e-02  &  9.18e-01  &  -2.52e-01  &  -1.20e-02  &  4.08e-02  &  -5.90e-04 \\
 1.09e-03  &  9.08e-04  &  -6.33e-04  &  -5.90e-04  &  7.92e-03  &  1.62e-01  &  9.18e-01  &  4.38e-02  &  2.21e-02  &  5.48e-02 \\
 3.79e-03  &  3.14e-03  &  -2.19e-03  &  -2.04e-03  &  2.74e-02  &  -9.34e-03  &  -1.06e-01  &  9.27e-01  &  1.35e-01  &  -1.21e-02 \\
 3.13e-03  &  2.59e-03  &  -1.81e-03  &  -1.69e-03  &  2.26e-02  &  3.53e-04  &  7.34e-03  &  -1.30e-01  &  9.56e-01  &  1.39e-01 \\
 4.47e-03  &  3.70e-03  &  -2.58e-03  &  -2.41e-03  &  3.23e-02  &  -1.05e-04  &  -5.21e-05  &  1.89e-04  &  -1.17e-04  &  1.00e+00 \\
\end{pmatrix},\\
B&=\begin{pmatrix}
 -4.58e-01  &  3.06e-17  &  3.30e-19 \\
 -1.69e-01  &  5.31e-02  &  -3.52e-05 \\
 2.77e-01  &  -4.41e-02  &  -1.41e-05 \\
 -2.98e-01  &  -3.67e-02  &  1.09e-05 \\
 1.97e-02  &  5.20e-01  &  3.09e-06 \\
 6.21e-04  &  -5.04e-03  &  4.69e-01 \\
 -2.82e-04  &  2.29e-03  &  7.61e-01 \\
 -9.76e-04  &  7.92e-03  &  3.26e-01 \\
 -8.06e-04  &  6.53e-03  &  -7.85e-02 \\
 -1.15e-03  &  9.33e-03  &  -1.56e-01 \\
\end{pmatrix},
\end{align*}}
results. 
The state and input constraints \eqref{eq:Constraints} read  
\begin{align*}
  -10\leq &\, x_i(t) \leq 10, \quad i= 1, \dots, 10,\\  
  -1\leq &\,u_j(t) \leq 1, \quad j= 1, \dots, 3,
\end{align*}
for this example. Furthermore, $Q = I^{n\times n}$, $R = 0.25I^{m\times m}$ and $N = 30$. The terminal weighting matrix $P$ is set to the solution of the discrete-time algebraic Riccati equation (DARE). The resulting quadratic program \eqref{eq:ParametricMPC} has $mN= 90$ decision variables and $q= 780$ inequality constraints. We reparametrize the inputs with the LQR controller proposed in \cite{Rossiter1998}. A condition number $\condNum{H}=2.51$ results from this reparametrization.
 
\paragraph*{{\bf MIMO75:}} MIMO75 differs from MIMO30 only with respect to the horizon length, which is set to $N = 75$ here.  
System matrices, constraints, weighting matrices and the input reparametrization are as in MIMO30. 
The resulting optimization problem \eqref{eq:ParametricMPC} has $mN= 225$ decision variables and $q= 1950$ inequality constraints.
 
\paragraph*{{\bf MIMORED30:}} 
System matrices, constraints, weighting matrices and the input reparametrization are the same as in MIMO30. In contrast to MIMO30, we here remove redundant constraints from \eqref{eq:ParametricMPC}.  
A constraint is called redundant if it never becomes active (see, for example, \cite[Sect.~4.1.1,~p.~128~ff.]{Boyd2009} or \cite[Def.~5,~p.~492~ff.]{Tondel2003}). We use the MPT toolbox~\cite{mpttoolbox} to identify redundant constraints.  
The resulting optimization problem \eqref{eq:ParametricMPC} has $mN= 90$ decision variables and $q= 556$ inequality constraints. 

\paragraph*{{\bf ACC25:}} This example is based on a discrete-time system that models an adaptive cruise control (ACC) \cite{Naus2010,Oliveri2011}. The ACC essentially controls the distance between a car equipped with ACC and the car in front of it. 
The two vehicles are referred to as the host and target, respectively. 
We implement the variant with a distance error penalty in the objective function~\cite{Oliveri2011}. The matrices of the system \eqref{eq:stateSpaceSystem} read 
\begin{equation} \label{eq:example_acc}
\begin{split}
A= \begin{pmatrix}1 & -T_s & 0 & 1.5T_s+\frac{1}{2}T_s^2 \\
0 & 1 & 0 & -T_s \\
0 & 0 & 1 & 0 \\
0 & 0 & 0 & 1
\end{pmatrix},\quad 
B=  \begin{pmatrix}0 \\0 \\0 \\ 1\end{pmatrix},
\end{split}
\end{equation}
where $T_s= 0.1\text{s}$ and  $x(t) = (e(t), v_r(t), v_t(t), a_h(t-1))^\prime$ 
with distance error $e(t)$, 
relative velocity $v_r(t) = v_t(t) - v_h(t)$, 
target vehicle velocity $v_t(t)$,
host vehicle velocity $v_h(t)$ and
host vehicle acceleration $a_h(t)$.  
The constraints \eqref{eq:Constraints} are as follows for ACC25:
\begin{equation}\label{eq:ACC_constraints}
  \begin{split}
\begin{array}{rcccl}
3.5+1.5(v_t-v_r)-200 &\leq& e(t) &\leq& 3.5+1.5(v_t-v_r),\\
v_t-50 &\leq& v_r(t) &\leq& v_t,\\
0 &\leq& v_t(t) &\leq& 50,\qquad \qquad\qquad \qquad \eqref{eq:ACC_constraints}\\
-3 &\leq& a_h(t) &\leq& 2, \\
-0.3 &\leq& u(t) &\leq& 0.3
\end{array}\\
  \end{split}\nonumber
\end{equation}
We briefly note that the state $v_t$ is not stabilizable. The model can be used to regulate the inter-vehicle distance and the relative speed, however~\cite{Naus2010,Oliveri2011}. 
The weighting matrices and the horizon are set to $Q=\text{diag}(2.5,5,0,1)$, $R=1$, $P = 0$ and $N = 25$, respectively.
The resulting quadratic program \eqref{eq:ParametricMPC} has $mN= 25$ decision variables and $q= 258$ constraints. Since the condition number of the matrix $H$ is $\condNum{H} = 4930.85$, we do not introduce an input reparametrization in this case.

\paragraph*{{\bf INPE50:}} We consider an inverted pendulum on a cart~\cite[p.~85~ff.]{Lunze2004}. The state vector reads $x(t) = \left(s(t),\varphi(t),\dot{s}(t),\dot{\varphi}(t)\right)^\prime$, where $s(t)$ is the position of the cart and $\varphi(t)$ is the pendulum angle. Zero-order hold discretization with $T_s=0.05\text{s}$ results in a system \eqref{eq:stateSpaceSystem} with matrices
\begin{align*}
A = 
\begin{pmatrix}
 1.00e+00  &  -1.07e-03  &  4.77e-02  &  -1.11e-05 \\
 0.00e+00  &  1.03e+00  &  4.73e-03  &  5.03e-02 \\
 0.00e+00  &  -4.22e-02  &  9.09e-01  &  -8.00e-04 \\
 0.00e+00  &  1.08e+00  &  1.87e-01  &  1.02e+00 \\
\end{pmatrix} ,\quad
B = \begin{pmatrix}
 3.63e-04 \\
 -7.53e-04 \\
 1.43e-02 \\
 -2.97e-02 \\
\end{pmatrix}. 
\end{align*}
The state and input constraints \eqref{eq:Constraints} read  
\begin{align*}
 \begin{split}
 \begin{array}{rcccl}
-1\leq &s(t)&\leq 1,\\
-\frac{\pi}{3} \leq &\varphi(t)& \leq \frac{\pi}{3},\\
-9\leq &\dot{x}(t)& \leq 9,\\
-2\pi\leq &\dot{\varphi}(t) &\leq 2\pi,\\
 -10\leq &u(t)&\leq 10.
\end{array}\\
  \end{split}\nonumber
\end{align*}
The weighting matrices are set to $Q = I^{n\times n}$, $R = 0.01 I^{m\times m}$, and $P$ is set to the result to the DARE. Choosing $N=50$ yields a quadratic program \eqref{eq:ParametricMPC} with $mN= 50$ decision variables and $q= 500$ inequality constraints. 
We reparametrize the inputs with the LQR controller proposed in~\cite{Rossiter1998}. A condition number $\condNum{H} = 1.00$ results after reparametrization.

\paragraph*{{\bf COMA40:}} This system models a linear chain of six masses connected to each other by springs, and to rigid walls on both ends of the chain~\cite{Wang2010,Wang2008}. All masses and all spring constants are set to unity. There exist three inputs $u_1, u_2, u_3$ that model forces between the first and second, third and fifth, and fourth and sixth mass, respectively. 
The resulting system has $12$ states and $3$ inputs. Discretizing with zero-order hold and sample time $T_s=0.5\text{s}$ yields a system of the form~\eqref{eq:stateSpaceSystem} with matrices

{\tiny
 \begin{align*}
A &= \left(\begin{array}{p{1.1cm}p{1.1cm}p{1.1cm}p{1.1cm}p{1.1cm}p{1.1cm}p{1.1cm}p{1.1cm}p{1.1cm}p{1.1cm}p{1.1cm}p{1.1cm}}
7.63e-01  &  1.15e-01  &  2.48e-03  &  2.09e-05  &  9.42e-08  &  2.63e-10  &  4.60e-01  &  1.98e-02  &  2.51e-04  &  1.51e-06  &  5.26e-09  &  1.20e-11 \\
 1.15e-01  &  7.65e-01  &  1.15e-01  &  2.48e-03  &  2.09e-05  &  9.42e-08  &  1.98e-02  &  4.60e-01  &  1.98e-02  &  2.51e-04  &  1.51e-06  &  5.26e-09 \\
 2.48e-03  &  1.15e-01  &  7.65e-01  &  1.15e-01  &  2.48e-03  &  2.09e-05  &  2.51e-04  &  1.98e-02  &  4.60e-01  &  1.98e-02  &  2.51e-04  &  1.51e-06 \\
 2.09e-05  &  2.48e-03  &  1.15e-01  &  7.65e-01  &  1.15e-01  &  2.48e-03  &  1.51e-06  &  2.51e-04  &  1.98e-02  &  4.60e-01  &  1.98e-02  &  2.51e-04 \\
 9.42e-08  &  2.09e-05  &  2.48e-03  &  1.15e-01  &  7.65e-01  &  1.15e-01  &  5.26e-09  &  1.51e-06  &  2.51e-04  &  1.98e-02  &  4.60e-01  &  1.98e-02 \\
 2.63e-10  &  9.42e-08  &  2.09e-05  &  2.48e-03  &  1.15e-01  &  7.63e-01  &  1.20e-11  &  5.26e-09  &  1.51e-06  &  2.51e-04  &  1.98e-02  &  4.60e-01 \\
 -8.99e-01  &  4.20e-01  &  1.93e-02  &  2.48e-04  &  1.50e-06  &  5.24e-09  &  7.63e-01  &  1.15e-01  &  2.48e-03  &  2.09e-05  &  9.42e-08  &  2.63e-10 \\
 4.20e-01  &  -8.80e-01  &  4.20e-01  &  1.93e-02  &  2.48e-04  &  1.50e-06  &  1.15e-01  &  7.65e-01  &  1.15e-01  &  2.48e-03  &  2.09e-05  &  9.42e-08 \\
 1.93e-02  &  4.20e-01  &  -8.80e-01  &  4.20e-01  &  1.93e-02  &  2.48e-04  &  2.48e-03  &  1.15e-01  &  7.65e-01  &  1.15e-01  &  2.48e-03  &  2.09e-05 \\
 2.48e-04  &  1.93e-02  &  4.20e-01  &  -8.80e-01  &  4.20e-01  &  1.93e-02  &  2.09e-05  &  2.48e-03  &  1.15e-01  &  7.65e-01  &  1.15e-01  &  2.48e-03 \\
 1.50e-06  &  2.48e-04  &  1.93e-02  &  4.20e-01  &  -8.80e-01  &  4.20e-01  &  9.42e-08  &  2.09e-05  &  2.48e-03  &  1.15e-01  &  7.65e-01  &  1.15e-01 \\
 5.24e-09  &  1.50e-06  &  2.48e-04  &  1.93e-02  &  4.20e-01  &  -8.99e-01  &  2.63e-10  &  9.42e-08  &  2.09e-05  &  2.48e-03  &  1.15e-01  &  7.63e-01 \\
\end{array}\right),\\
B &= \left(\begin{array}{p{1.1cm}p{1.1cm}p{1.1cm}}
 1.17e-01  &  2.11e-05  &  9.48e-08 \\
 -1.17e-01  &  2.52e-03  &  2.11e-05 \\
 -2.50e-03  &  1.20e-01  &  2.52e-03 \\
 -2.10e-05  &  5.02e-13  &  1.20e-01 \\
 -9.45e-08  &  -1.20e-01  &  9.48e-08 \\
 -2.64e-10  &  -2.52e-03  &  -1.20e-01 \\
 4.40e-01  &  2.51e-04  &  1.51e-06 \\
 -4.40e-01  &  1.98e-02  &  2.51e-04 \\
 -1.96e-02  &  4.60e-01  &  1.98e-02 \\
 -2.50e-04  &  1.20e-11  &  4.60e-01 \\
 -1.50e-06  &  -4.60e-01  &  1.51e-06 \\
 -5.25e-09  &  -1.98e-02  &  -4.59e-01 \\
\end{array}\right).
\end{align*}
}
The state and input constraints~\eqref{eq:Constraints} read 
\begin{align*}
  -4\leq &x_i(t)\leq 4, \quad i= 1, \dots, 12\\
  -0.5 \leq &u_j(t)\leq 0.5, \quad j= 1, \dots, 3,
\end{align*}
respectively. We choose $Q = I^{n\times n}$, $R = I^{m\times m}$ and set $P$ to the solution of the DARE. 
The resulting quadratic program~\eqref{eq:ParametricMPC} has $mN= 120$ decision variables and $q=1200$ inequality constraints for a horizon of $N= 40$. After reparametrization with the LQR controller proposed in \cite{Rossiter1998}, a condition number of $\condNum{H} = 1.47$ results.

\subsection{Brief description of the solvers}\label{subsec:solvers}
A brief comparison of the solvers is given in Tab.~\ref{tab:SummaryOptAlgs}. Additional comments are given below.    
\begin{table*}[b]
\centering
\caption{Summary of the optimization algorithms.}~\label{tab:SummaryOptAlgs}
\begin{tabular}{l||l|c|c}
  {\bf Name} & {\bf Type} & {\bf Version} & {\bf References} \\\hline \hline
  int-pnt-cvx & primal-dual interior-point algorithm& 2013a & \cite{OptToolBox}\\\hline 
  act-set & primal active-set algorithm &R2013a& \cite{OptToolBox} \\ \hline
  qpip & primal-dual interior-point algorithm &2.0&\cite{qpc} \\ \hline
  qpas & dual active-set algorithm&2.0&\cite{qpc} \\\hline
  OOQP & primal-dual interior-point algorithm &0.99.24& \cite{OOQP,Gertz2003} \\\hline 
  MOSEK & self-dual interior-point algorithm with presolve phase &  7 Build 121 &\cite{MOSEK} 
\end{tabular}
\end{table*}

\paragraph{{\bf Matlab interior point and active-set solvers~\cite{OptToolBox}:}} The Matlab Optimization Toolbox provides an interior point solver for convex QPs and an active set solver for QPs, which we refer to by int-pnt-cvx and act-set for short. We select these solvers, because they are easy to use and widely available. Consequently, our results obtained with these solvers can be verified particularly easily. 

\paragraph{{\bf qpip and qpas~\cite{qpc}:}} The solvers qpip and qpas from the QPC library~\cite{qpc} implement a primal-dual interior-point and a dual active-set algorithm, respectively, for strictly convex QPs. They are implemented in C and can be called from Matlab. We include these solvers, because they are claimed to be several orders of magnitude faster than the solvers from the Matlab Optimization Toolbox~\cite{qpc}. Moreover, qpas is included, because one of the reviewers of~\cite{Jost2014} suggested to apply constraint removal to a dual active-set method. Note that these solvers are only available as binaries. Since constraint removal does not require to change the solver itself, this is not a restriction, however.

\paragraph{{\bf OOQP~\cite{OOQP,Gertz2003}:}} OOQP implements a primal-dual interior-point algorithm for convex QPs. 
Specifically, we use the option that is based on Mehrotra’s predictor-corrector algorithm with Gondzio’s multiple corrections (see \cite{OOQP} and the references therein).  
We include OOQP in order to test constraint removal on an implementation that does not involve Matlab.\footnote{We note that OOQP is an interesting option, because it is documented very well and very flexible w.r.t.\ the underlying linear algebra and compilation for gcc-supporting platforms beyond desktop PCs.}

\paragraph{{\bf MOSEK~\cite{MOSEK}:}} Among other solvers, MOSEK provides a self-dual interior-point solver for convex QPs. 
MOSEK's presolve phase presumably removes redundant constraints \cite{MOSEK}, but no algorithmic details are given. We select this solver to demonstrate that constraint removal can be used with commercial solvers. We claim MOSEK is a typical commercial solver in that it has been used to solve a number of academic and industrial problems (see the references in \cite{MOSEK}) on the one hand, and little information on the implemented algorithms are available on the other hand. The results presented in Sect.~\ref{sec:results} show that the lack of information on the algorithmic details of a solver does not impede constraint removal at all. Finally, we note that we call MOSEK from its Matlab interface.

\section{Results}\label{sec:results}
We compare computational times using the cumulative distribution functions (cdf) $h_{\text{cdf}}(t_{\text{MPC}})$ in Sect.~\ref{subsec:CDFs}. The cumulative distribution function $h_{\text{cdf}}(t)$ is defined as the fraction of QPs in which the control law is found in time $t$ or less. For each of the 36 combinations of the six examples and six solvers, we determine and compare the cumulative distribution function that results with and without constraint removal. 

A higher level summary and an interpretation of the results are given in Sects.~\ref{subsec:SummaryComputationTimes} and \ref{subsec:Interpretation}, respectively.

\subsection{Cumulative distribution functions of the computational times}\label{subsec:CDFs}
The cdfs for all 36 cases are shown in Figs.~\ref{fig:resMIMO}--\ref{fig:resCOMA}. The intervals on the abscissae of Figs.~\ref{fig:resMIMO}--\ref{fig:resCOMA} are chosen such that at least the range $[0, 0.99]$ of the cdf is displayed. Enlargements are added where appropriate. 
Note that the time $t_{\text{MPC}}$ such that $h_{\text{cdf}}(t_{\text{MPC}}) = 1$ is the maximal computational time required in that case. 

\begin{figure}[H]
  \centering	
  \includegraphics[width=\textwidth]{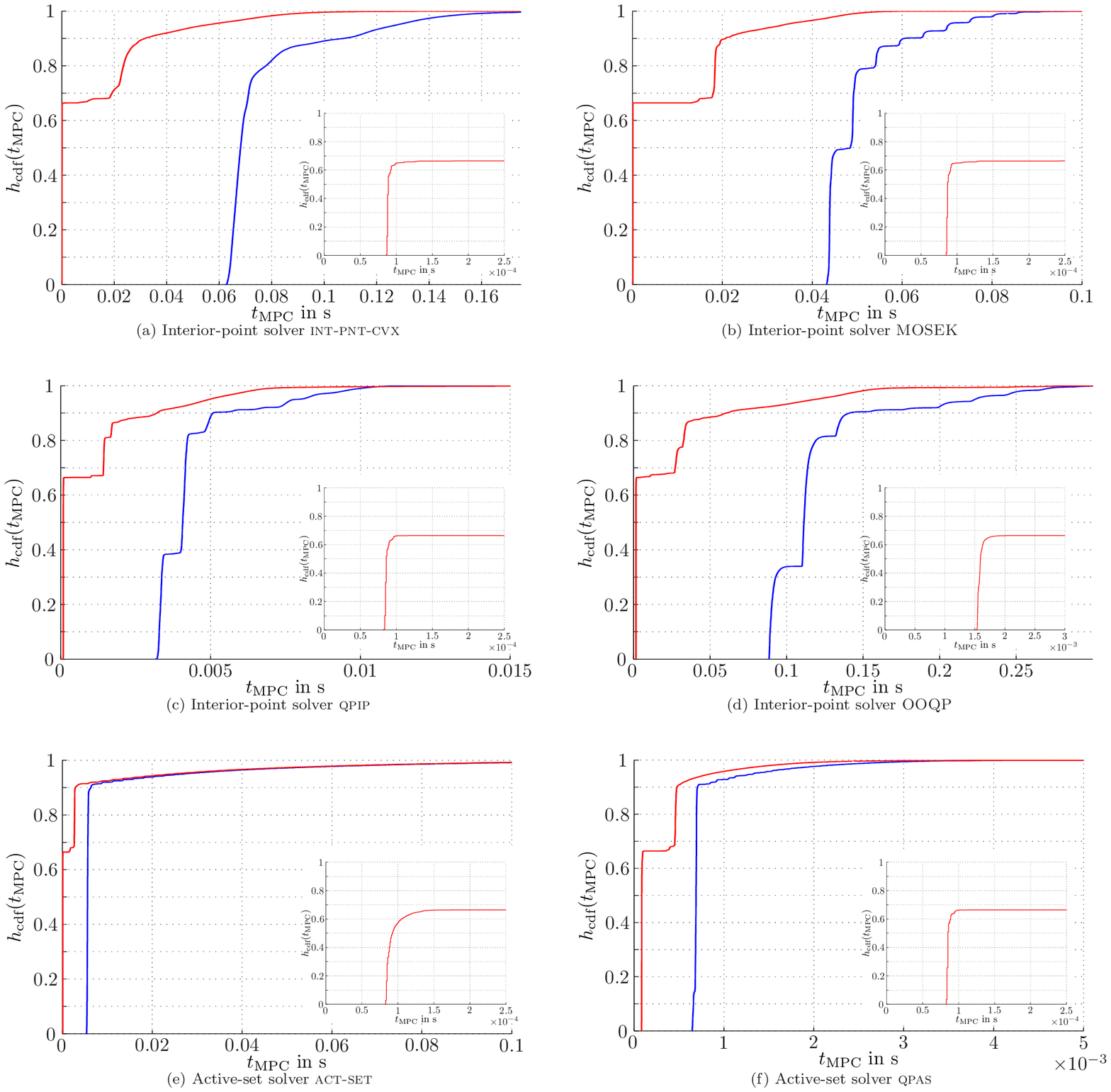}
  \caption{Cumulative distribution functions for the MIMO30 example (full-MPC in blue, CR-MPC in red).}\label{fig:resMIMO}
\end{figure}

\begin{figure}
  \centering	
  \includegraphics[width=\textwidth]{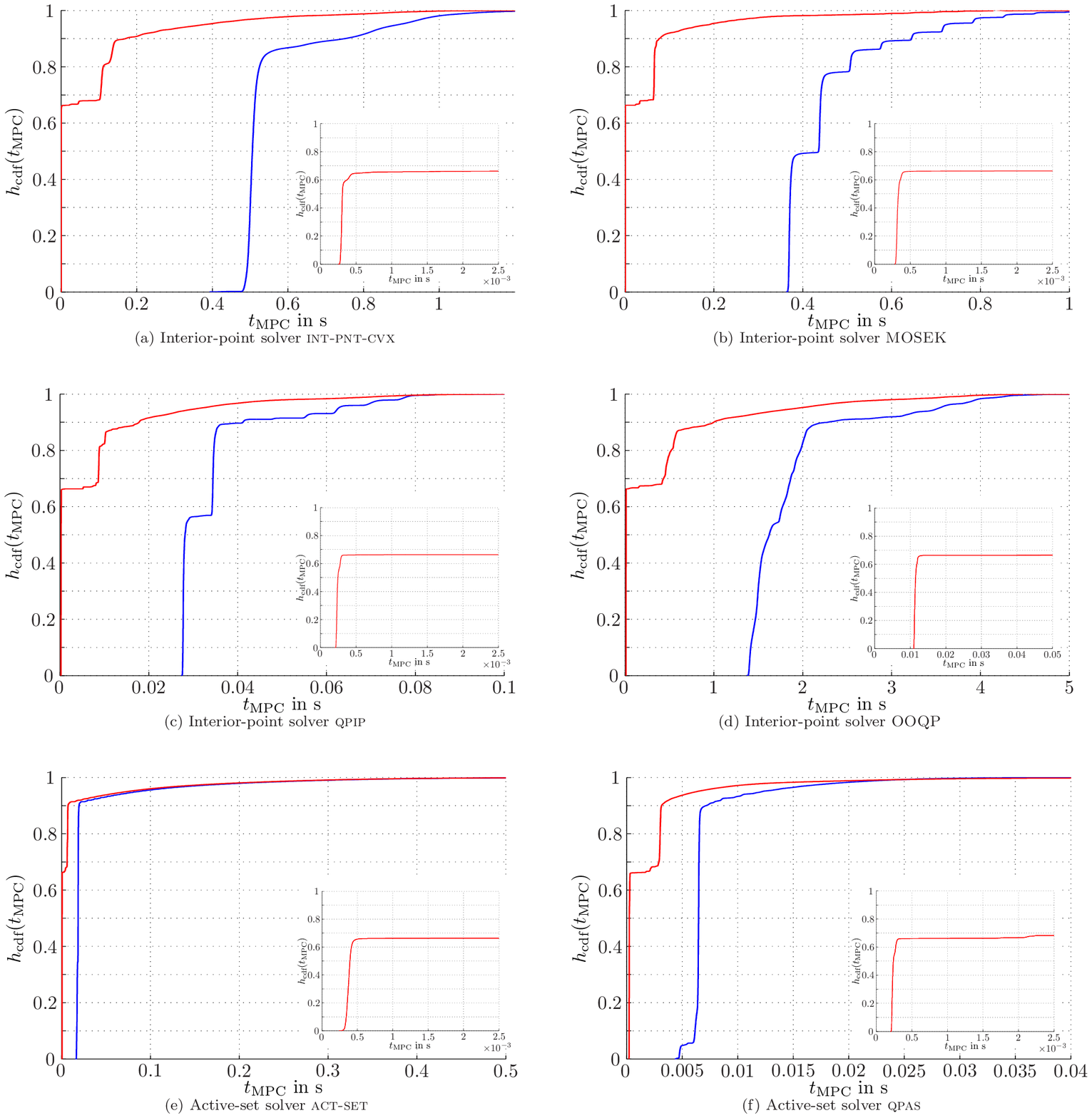}
  \caption{Cumulative distribution functions for the MIMO75 example (full-MPC in blue, CR-MPC in red).}\label{fig:resMIMO75}
\end{figure}

\subsubsection*{Results for MIMO30, Fig.~\ref{fig:resMIMO}}
  \begin{itemize}
  \item Approximately 65\% of the QPs are detected to be unconstrained by CR-MPC. Consequently, no optimization problem is solved at all in these cases by CR-MPC, but the optimal control law of the unconstrained case is evaluated immediately (line 4 of Algorithm~\ref{alg:constraintRemovalAlgorithm}). CR-MPC therefore provides the optimal input sequence very quickly, which results in the leftmost shoulders on all red curves in Figs.~\ref{fig:resMIMO}a--f. Moreover, it is evident from Figs.~\ref{fig:resMIMO}a--c, e, f (but not d) that the computation time does not depend on the QP solver if the QP is detected to be unconstrained. The data shown in Figs.~\ref{fig:resMIMO}a--c, e, f have all been obtained in Matlab and therefore all result in the same leftmost shoulder on the red curves. The case shown in Figs.~\ref{fig:resMIMO}d is implemented in C/C++ and therefore results in a leftmost edge that also rises to about 65\%, but at a different computational time.
  
  \item For both the interior-point solvers (Figs.~\ref{fig:resMIMO}a-d) and the active-set solvers (Figs.~\ref{fig:resMIMO}e-f), MPC with constraint removal outperforms MPC without constraint removal in the sense that the cdf for CR-MPC always lies to the left and above that for full-MPC. 
      
  \item The maximal computation time required by CR-MPC is smaller than that required by full-MPC for all QP solvers. The difference is pronounced for the interior-point solvers but small for the active-set solvers. 
\end{itemize}

\subsubsection*{Results for MIMO75, Fig.~cdf's\ref{fig:resMIMO75}}
Recall this example differs from MIMO30 in that the horizon is longer here ($N= 75$). Consequently, the number of decision variables and constraints is larger in MIMO75 than in MIMO30 (by a factor of approximately 2.5, see Tab.~\ref{tab:systems}). 
\begin{itemize}
  \item Approximately 65\% of the QPs are detected to be unconstrained by CR-MPC. Consequently, no optimization problem needs to be solved at all. See the first comment on the MIMO30 results for a more detailed discussion. 
  \item MPC with constraint removal outperforms MPC without constraint removal for all solvers in the sense stated in the second comment on the MIMO30 results. 
  \item All curves have approximately the same shape as the corresponding ones for MIMO30, but they are shifted to larger values of $t_{\text{MPC}}$ here. This result is consistent with the fact that both MIMO30 and MIMO75 are based on the same system \eqref{eq:stateSpaceSystem} and the same constraints \eqref{eq:Constraints}, but the number of constraints is larger here due to the larger horizon. 
  \item The maximal computation times required by CR-MPC are smaller than for full-MPC. The same observations hold as stated in the third comment on the MIMO30.
\end{itemize}

\begin{figure}
  \centering	
  \includegraphics[width=\textwidth]{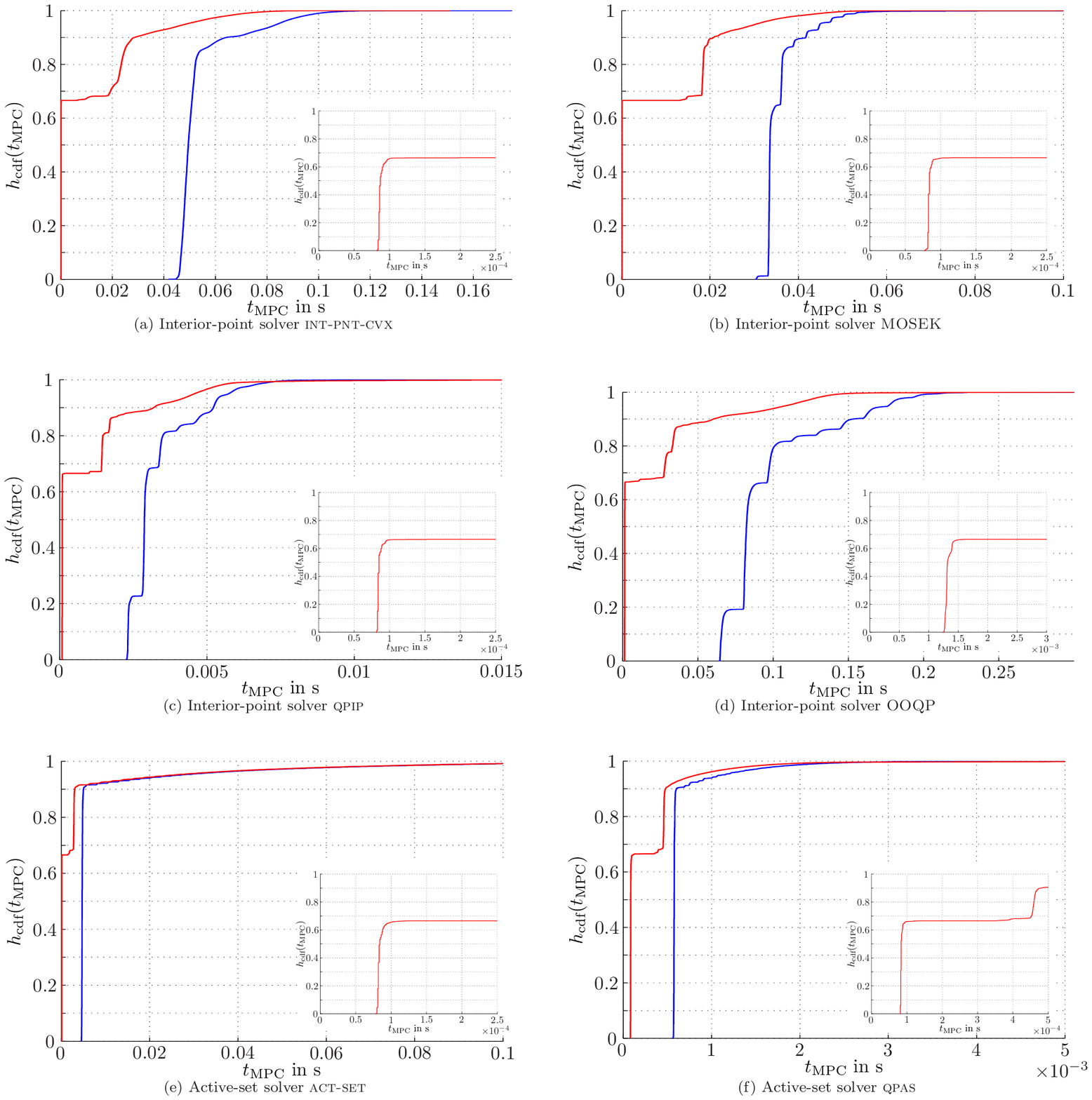}
  \caption{Cumulative distribution functions for the MIMORED30 example (full-MPC in blue, CR-MPC in red).}\label{fig:resMIMORED}
\end{figure}

\subsubsection*{Results for MIMORED30, Fig.~\ref{fig:resMIMORED}}
Recall MIMORED30 differs from MIMO30 only in that the redundant constraints have been removed in MIMORED30.
\begin{itemize}
 \item Approximately 65\% of the optimization problems are detected to be unconstrained by CR-MPC. 
 Consequently, CR-MPC is considerably faster than full-MPC, since no QP is solved at all in CR-MPC.
 See the first comment on the MIMO30 results for a more detailed discussion. 
 
 \item MPC with constraint removal outperforms MPC without constraint removal for all solvers in the sense stated in the second comment on the MIMO30 results. One of the reviewers of \cite{Jost2014} conjectured constraint removal would have a weak impact, if any, on problems in which redundant constraints are removed a priori. Note that this is clearly not true in this example. In other words, constraint removal does not just remove redundant constraints. 
 
 \item We can analyze the role of redundant constraints in more detail by comparing the results of MIMORED30 to those of MIMO30.
  Consider the cdfs of the full-MPC cases (blue curves in Fig.~\ref{fig:resMIMORED}a--f) first. These cdfs have approximately the same shape as the corresponding ones for MIMO30 (blue curves in Fig.~\ref{fig:resMIMO}a--f), but they are shifted to smaller values of $t_{\text{MPC}}$ here. 
  These shifts result, since the redundant constraints have been removed here. Table~\ref{tab:shift} shows some more precise data on the shifts. Specifically, Table~\ref{tab:shift} lists the times $t_{\text{MPC}}$ for which the cdfs attain the value $0.7$ (data for full-MPC and MIMO30 in column 2, data for full-MPC and MIMORED30 in column 3).\footnote{The value $0.7$ is arbitrary. A value larger than $0.65$ should be chosen, because about 65\% of the QPs are detected to be unconstrained by CR-MPC.}
 The fourth column lists how large the shift of the cdf is. Observe the figures range from about 15\% to about 29\% if constraint removal is not applied. 
 \newline
 Now consider the corresponding diagrams (red curves in Fig.~\ref{fig:resMIMORED}a-f) and data (columns 5-7 in Table~\ref{tab:shift}) for CR-MPC, i.e.\ for the cases with constraint removal. These cdfs also have approximately the same shape as the corresponding cdfs resulting for the MIMO30 example (red curves in Fig.~\ref{fig:resMIMO}a--f), but the shift is much smaller than for full-MPC. The percentages range from about 0.1\% to about 4\% (column 7 of Table~\ref{tab:shift}) if constraint removal is applied compared to 15\% to 29\% found above for full-MPC.
 \newline
 Essentially, this comparison indicates that removing redundant constraints has a strong impact on the computation times of full-MPC, but only a much smaller impact on those of CR-MPC. In other words, redundant constraints need not be removed if CR-MPC is used, since CR-MPC removes many of the redundant constraints before invoking the QP solver anyhow. 
\begin{table}
\centering
  \setlength{\tabcolsep}{0.35mm}
 \scriptsize
 \caption{Calculation times $t_{\text{MPC}}$ resulting at $h_{\text{cdf}}(t_{\text{MPC}})=0.7$. All absolute values are given in seconds.}\label{tab:shift}
\begin{tabular}{c||c|c|c||c|c|c}
Solver & 
(full-MPC, MIMO30) &
(full-MPC, MIMORED30) &
shift &
(CR-MPC, MIMO30) &
(CR-MPC, MIMORED30) &
shift 
\\\hline \hline
int-pnt-cvx &
$71.65\cdot 10^{-3}$&
$51.31\cdot 10^{-3}$&
$-28.39\%$&
$19.78\cdot 10^{-3}$&
$19.03\cdot 10^{-3}$&
$-3.79\%$
\\ \hline
MOSEK &
$50.72\cdot 10^{-3}$&
$36.11\cdot 10^{-3}$&
$-28.81\%$&
$18.23\cdot 10^{-3}$&
$18.19\cdot 10^{-3}$&
$-0.22\%$
\\\hline
qpip &
$4.161\cdot 10^{-3}$&
$3.373\cdot 10^{-3}$&
$-18.94\%$&
$1.428\cdot 10^{-3}$&
$1.426\cdot 10^{-3}$&
$-0.14\%$
\\\hline
OOQP &
$113.8\cdot 10^{-3}$&
$96.95\cdot 10^{-3}$&
$-14.81\%$&
$27.48\cdot 10^{-3}$&
$27.21\cdot 10^{-3}$&
$-0.99\%$
\\\hline
act-set &
$5.937\cdot 10^{-3}$&
$4.782\cdot 10^{-3}$&
$-19.45\%$&
$2.870\cdot 10^{-3}$&
$2.851\cdot 10^{-3}$&
$-0.69\%$
\\\hline
qpas &
$0.693\cdot 10^{-3}$&
$0.579\cdot 10^{-3}$&
$-16.45\%$&
$0.461\cdot 10^{-3}$&
$0.459\cdot 10^{-3}$&
$-0.22\%$
\\
\end{tabular}
\end{table} 
\end{itemize}

\subsubsection*{Results for ACC25, Fig.~\ref{fig:resACC}} 
\begin{itemize}
  \item Approximately 55\% of the QPs are detected to be unconstrained by CR-MPC. 
  Consequently, CR-MPC is considerably faster than full-MPC, since no QP is solved at all in CR-MPC.
  See the first comment on the MIMO30 results for a more detailed discussion. 

  \item For all interior-point solvers and the active-set solver act-set CR-MPC outperforms full-MPC in the sense stated in the second comment on the MIMO30 results. For the dual active-set solver qpas this claim does not hold in this example. See the next comment. 
  
  \item Consider the cdfs for the qpas solver shown in Fig.~\ref{fig:resACC}f in more detail. For calculation times in the range $0.18\cdot10^{-3}\text{s}<t_{\text{MPC}}<0.3\cdot10^{-3}\text{s}$ the cdf of full-MPC is larger than that for CR-MPC while the converse holds for $t_{\text{MPC}}< 0.18\cdot10^{-3}\text{s}$ and $t_{\text{MPC}}>0.3\cdot10^{-3}\text{s}$.
  This implies there exist QPs for which the additional calculations required to generate the reduced QP (specifically, by lines 2 and 6 in Alg.~\ref{alg:constraintRemovalAlgorithm}) require more time than saved due to the reductions.  We claim this happens for two reasons. Most importantly, the speed-up due constraint removal must be expected to be much smaller for active-set methods than for interior-point solvers. This is, roughly speaking, due to the fact that interior-point solvers always operate on all constraints, while active-set solvers operate on a subset thereof (see Sect.~\ref{subsec:Interpretation} for more detailed comments). Secondly, the smaller the number of constraints in \eqref{eq:ParametricMPC}, the smaller the impact of constraint removal. In fact, ACC25 is the smallest example treated here (see Tab.~\ref{tab:systems}). 

  \item The maximal computation time required by CR-MPC is smaller than for full-MPC. See the third comment on the MIMO30 example.
\end{itemize}

\subsubsection*{Results for INPE50, Fig.~\ref{fig:resINPE}}
  \begin{itemize}
  \item Approximately 88\% of the QPs are detected to be unconstrained by CR-MPC. Consequently, no optimization problem needs to be solved at all. See the first comment on the MIMO30 results for a more detailed discussion. 
  
  \item Again, CR-MPC outperforms full-MPC for all solvers in the sense that the resulting cdfs lie to the left and above of the cdfs resulting for full-MPC (cf. Fig.~\ref{fig:resINPE}a-f).
  
  \item The maximal computation time required by CR-MPC is smaller than for full-MPC. The same observation holds as stated in the third comment on the MIMO30 example.
\end{itemize}
 
\paragraph{Results for COMA40, Fig.~\ref{fig:resCOMA}} 
\begin{itemize}
  \item Approximately 65\% of the QPs are detected to be unconstrained by CR-MPC. Consequently, no optimization problem needs to be solved at all. See the first comment on the MIMO30 results for a more detailed discussion. 
  
  \item CR-MPC outperforms full-MPC for all solvers (cf. Fig.~\ref{fig:resCOMA}a-f) in the sense discussed in the second comment on the MIMO30 results.
  
  \item The maximal computation time required by CR-MPC is smaller than for full-MPC. The same observation holds as stated in the third comment on the MIMO30 example.   
\end{itemize}

\begin{table}[t]
 \setlength{\tabcolsep}{0.35mm}
 \scriptsize
\centering
\caption{Difference between the average computation times resulting for CR-MPC and full-MPC for all example-solver combinations. Absolute values are given in milliseconds.}~\label{tab:crAverageComputationTime}
\begin{tabular}{c||c c|c c|c c|c c|c c|c c||c}
 Problem &  \multicolumn{2}{c|}{int-pnt-cvx} & \multicolumn{2}{c|}{MOSEK} & \multicolumn{2}{c|}{qpip} & \multicolumn{2}{c|}{OOQP} &\multicolumn{2}{c|}{act-set} & \multicolumn{2}{c||}{qpas} & average \\ \hline \hline
  MIMO30 &  -64.59 &(  -85.10\%) &   -44.29 &(  -84.15\%) &    -3.29 &(  -76.95\%) &   -97.12 &(  -82.28\%) &    -4.86 &(  -51.84\%) &    -0.50 &(  -64.19\%)  & -74.09\% \\\hline
  MIMO75 & -481.11 &(  -87.06\%) &  -413.09 &(  -90.90\%) &   -27.51 &(  -81.60\%) & -1493.83 &(  -81.29\%) &   -16.15 &(  -54.33\%) &    -5.27 &(  -74.10\%) &-78.21\% \\\hline
  MIMORED30 &   -42.57 &(  -80.61\%) &   -28.18 &(  -78.82\%) &    -2.35& (-71.45\%) &   -76.14 &(  -79.83\%) &    -3.71 &(  -45.18\%) &    -0.38 &(  -57.16\%)  & -68.46\% \\\hline
  ACC25  &   -26.04 &(  -74.84\%) &    -6.50 &(  -54.74\%) &    -0.80 &(  -69.58\%) &    -7.32 &(  -78.23\%) &    -1.49 &(  -49.55\%) &    -0.05 &(  -23.92\%) &-58.48\% \\\hline
  INPE50 &   -27.97& (  -95.47\%) &   -29.41& (  -93.77\%) &    -1.78& (  -89.32\%) &   -22.69& (  -94.61\%) &    -3.55& (  -90.02\%) &    -0.21& (  -66.45\%) &-88.27\% \\\hline
  COMA40 &  -163.34 &(  -82.76\%) &  -128.73 &(  -82.86\%) &    -9.36 &(  -75.32\%) &  -304.98 &(  -85.74\%) &    -7.81 &(  -16.45\%) &    -1.32 &(  -55.11\%) &-66.37\%\\  \hline \hline
   average&  -134.27 &( -84.31\%)  & -108.37 & (-80.87\%) &  -7.53 & (-76.98\%)& -333.68 & (-83.66\%) & -6.26 &  (-51.23\%) &  -1.29 & (-56.82\%)\\
\end{tabular}
\end{table}
\subsection{Analysis of the computation times}\label{subsec:SummaryComputationTimes}
\begin{table}[b]
\centering
\caption{Fraction of QPs which have been solved faster with CR-MPC than the {\it fastest} time required to find the control law with full-MPC.}~\label{tab:crFasterThanFastestFull}
\begin{tabular}{c||c|c|c|c|c|c}
Problem &  int-pnt-cvx & MOSEK & qpip & OOQP & act-set & qpas  \\ \hline \hline
 MIMO30    & 95.90 \% &    97.02 \% &    89.66 \% &    92.30 \% &    91.47 \% &    93.07 \% \\ \hline
 MIMO75    & 95.24 \% &    98.05 \% &    93.94 \% &    92.44 \% &    91.57 \% &    92.94  \% \\ \hline
 MIMORED30 & 93.56 \% &    94.91 \% &    87.98 \% &    90.60 \% &    91.57 \% &    91.94 \% \\ \hline
 ACC25     & 75.86 \% &    63.27 \% &    94.07 \% &    97.80 \% &    92.70 \% &    53.75 \% \\ \hline
 INPE50    & 99.30 \% &    99.29 \% &    96.27 \% &    99.64 \% &    98.60 \% &    98.51 \% \\ \hline
 COMA40    & 92.59 \% &    93.35 \% &    89.47 \% &    97.65 \% &    77.09 \% &    79.14 \% 
\end{tabular}
\end{table}
Table~\ref{tab:crAverageComputationTime} shows the difference between the average computation times resulting for CR-MPC and full-MPC for each of the 36 example-solver combinations. Consider the results for the combination (MIMO30, int-pnt-cvx), for example (top left entry in Tab.~\ref{tab:crAverageComputationTime}). This figure means that the QPs are solved faster with CR-MPC than with full-MPC by 64.59 milliseconds on average. This equals a relative reduction of 85.10\%. 

Consider the second column of Tab.~\ref{tab:crAverageComputationTime}, which corresponds to the int-pnt-cvx solver. Here, the relative reduction of the computation time varies between 75\% and 95\% for all examples. Similar reductions result for the interior-point solvers qpip (col. 4, 69\%--89\%) and OOQP (col.\ 5, 78\%--95\%). For the interior-point solver MOSEK (col.\ 3) the minimal reduction (55\%) is slightly smaller, while the maximal reduction (94\%) is of the same order as with OOQP and qpip. For the active-set solvers act-set (col.\ 6) and qpas (col.\ 7) the reduction varies between 18\% and 90\% respectively 24\% and 75\%. Note that the minimal reduction is smaller for the active-set solvers than for the interior-point solvers. 
The last column of Tab.~\ref{tab:crAverageComputationTime} lists the average reduction for each example, which varies between 58\% and 88\%. The smallest figure results for the ACC25 example. This system is the smallest with respect to the number of constraints and decision variables among all examples (cf. Tab.~\ref{tab:systems}).

The last row of Tab.~\ref{tab:crAverageComputationTime} shows the average reduction for each solver. This figure varies between 78\% and 84\% for the interior-point solvers, and between 51\% and 57\% for the active-set solvers. 
These results confirm that constraint removal accelerates interior-point solvers more strongly than active-set solvers (see the discusssion in Sect.~\ref{subsec:Interpretation}).

Table~\ref{tab:crFasterThanFastestFull} lists the fractions of QPs that are solved in less time with CR-MPC than the \textit{shortest time} required by full-MPC \textit{among all QPs}. 
For example, the figure 95.9\% in the top left corner means that in about 96\% of all QPs the int-pnt-cvx solver with constraint removal requires less computational time than the shortest time required among all QPs with the same solver but without constraint removal. For completeness we note that this value can also be found in Fig.~\ref{fig:resMIMO}a: Consider a line parallel to the ordinate which cuts the abscissa at the smallest computation time achieved by full-MPC (blue curve). The entry listed in Tab.~\ref{tab:crFasterThanFastestFull} is equal to the value of the cdf resulting for CR-MPC (red curve) at the intersection with this line. 

The smallest entries in Tab.~\ref{tab:crFasterThanFastestFull} occur for the smallest system (ACC25). The three smallest figures arise for the combinations (ACC25, int-pnt-cvx), (ACC25, MOSEK) and (ACC25, qpas), which amount to about 75\%, 63\% and 54\%, respectively. Even in the worst case, CR-MPC is faster in 54\% than the fastest full-MPC run among all runs of the corresponding case.  Note that in 30 out of the 36 combinations listed in Tab.~\ref{tab:crFasterThanFastestFull} this figure is even larger than 90\%. 

\subsection{Brief interpretation of the results}\label{subsec:Interpretation}
Larger reductions in computational time result for interior-point solvers than for active-set solvers. We give an informal explanation of this result in this section.  

Active-set solvers search for the optimal active set by generating candidate active sets, which generally are subsets of the set of all constraints. A candidate optimal point is found by solving a quadratic subproblem, in which the constraints in the current candidate active set are treated as equality constraints and the remaining constraints are ignored. The candidate optimal point is then modified such that feasibility with respect to the remaining constraints is maintained. (Primal active-set solvers generate candidate optimal points that are feasible with respect to the primal problem~\cite[Chpt.~16.5,~p.~467~ff.]{Nocedal2000}, while dual active-set solvers generate those that are feasible with respect to the dual problem~\cite{Goldfarb1983}.) 
If constraints are removed by constraint removal, feasibility has to be enforced with respect to a smaller number of constraints only. The computational cost of maintaining feasibility is therefore reduced by constraint removal. Note that this effect remains even if the sequence of intermediate QPs solved by the active-set solver is not affected by constraint removal, because the removed constraints are never considered in any candidate active set.

Primal and dual active-set methods generate initial candidate active sets differently. This step is known to be time-consuming for primal active-set solvers. It may require up to 50\% of the time taken to solve the quadratic program~\cite{Goldfarb1983}. 
Dual active-set solvers, in contrast, can start from the solution of the unconstrained optimization problem, which is guaranteed to be dual feasible. We conjecture this difference is one of the reasons why constraint removal has a stronger effect on primal than on dual active-set solvers. 
It should be remarked that dual active-set solvers in particular may temporarily add constraints to the candidate set that are inactive at the optimal solution. Clearly, such constraints are later removed from the candidate set before the algorithm converges. Therefore, constraint removal results in a reduction even for dual active-set algorithms because it reduces the overall number of constraints and prevents temporary inclusion in the candidate set of (some of the) constraints that will be inactive at the optimal solution.

In contrast to active-set solvers, interior-point solvers do not operate on subsets of the constraints but on all constraints, because inequality constraints are rewritten as equality constraints with slack variables. 
The main effort during an iteration of an interior-point solver is to solve a linear system of equations that involves all inequality constraints~\cite[Chpt.~11.8,~p.~615~ff.]{Boyd2009}. These linear systems are obviously simplified if the number of constraints is reduced a-priori by constraint removal.

\section{Conclusion}\label{sec:conclusion}
Constraint removal accelerates MPC as anticipated. 
In 34 out of 36 example-solver combinations the average computation time is reduced considerably (45\%--95\%). The reduction is smaller (17\% and 23\%) in the remaining two combinations, but still significant (see Tab.~\ref{tab:crAverageComputationTime}). 
Larger reductions result for the interior-point solvers than for the active-set solvers. A reduction can be achieved even for the dual active-set solver, which was expected to be affected the least for algorithmic reasons (see Sect.~\ref{subsec:Interpretation}). 

\begin{figure}
  \centering	
  \includegraphics[width=\textwidth]{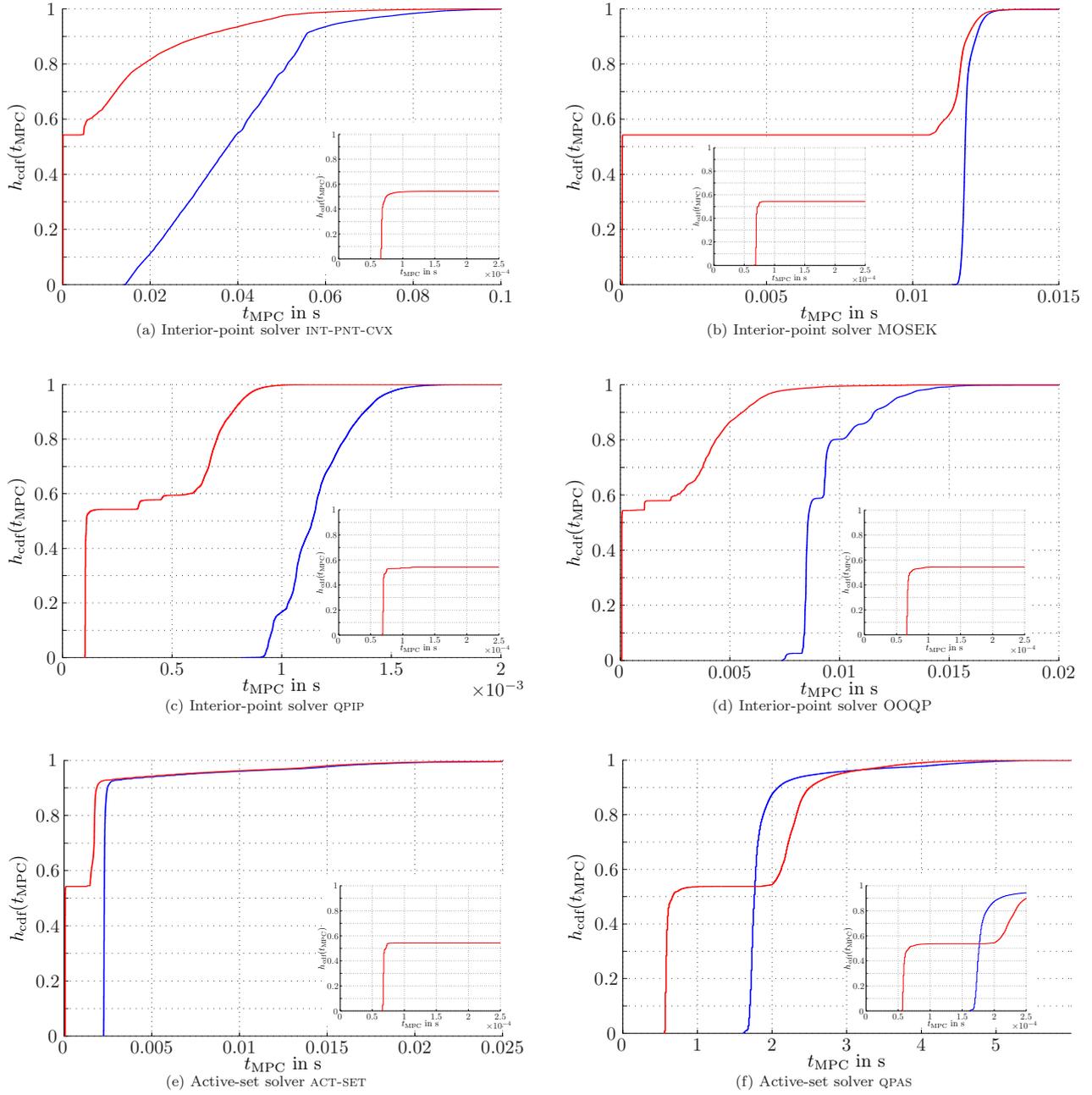}
  \caption{Cumulative distribution functions for the ACC25 example (full-MPC in blue, CR-MPC in red).}\label{fig:resACC}
\end{figure}

\begin{figure}
  \centering	
  \includegraphics[width=\textwidth]{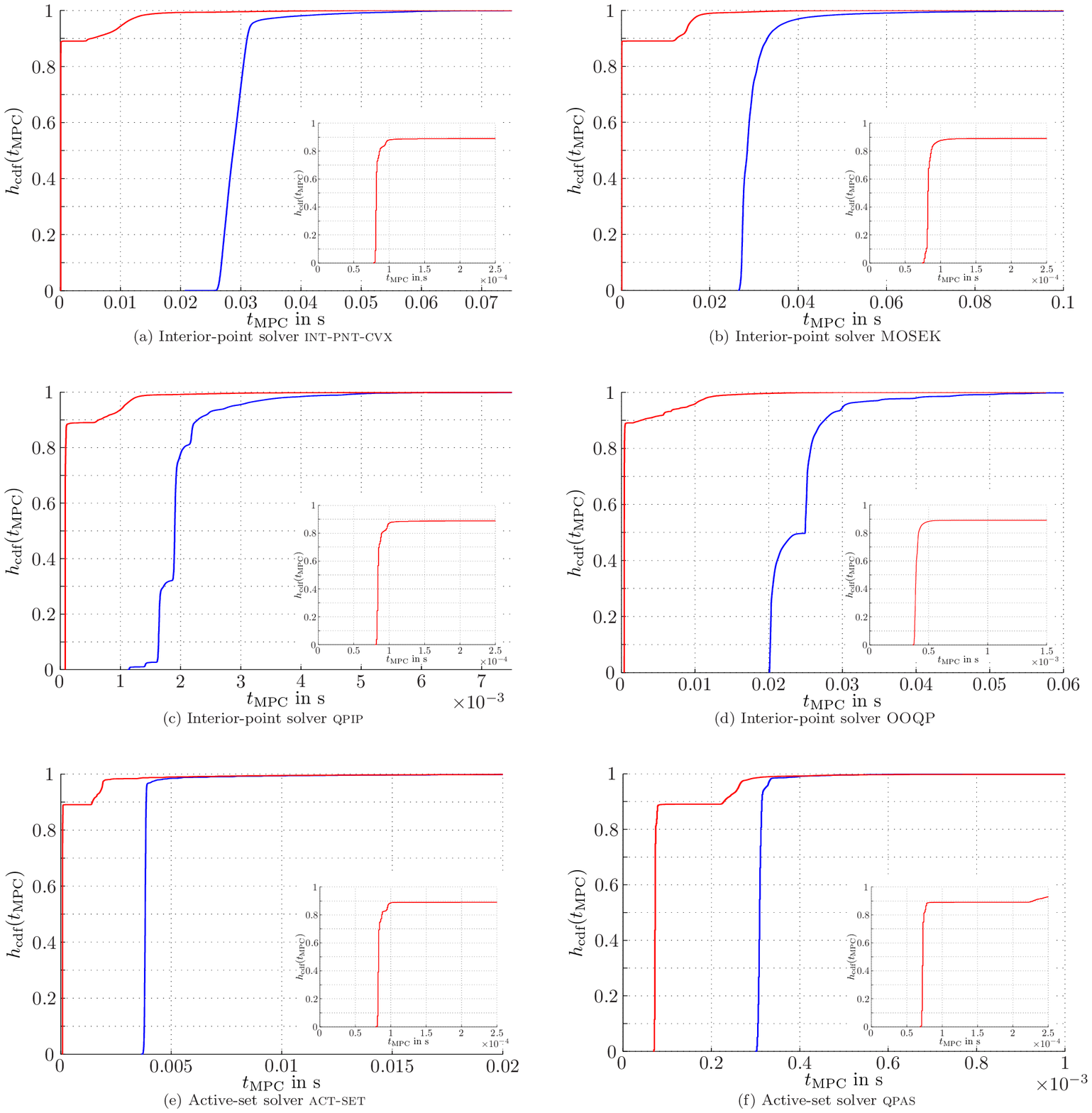}
  \caption{Cumulative distribution functions for the INPE50 example (full-MPC in blue, CR-MPC in red).}\label{fig:resINPE}
\end{figure}

\begin{figure}
  \centering	
  \includegraphics[width=\textwidth]{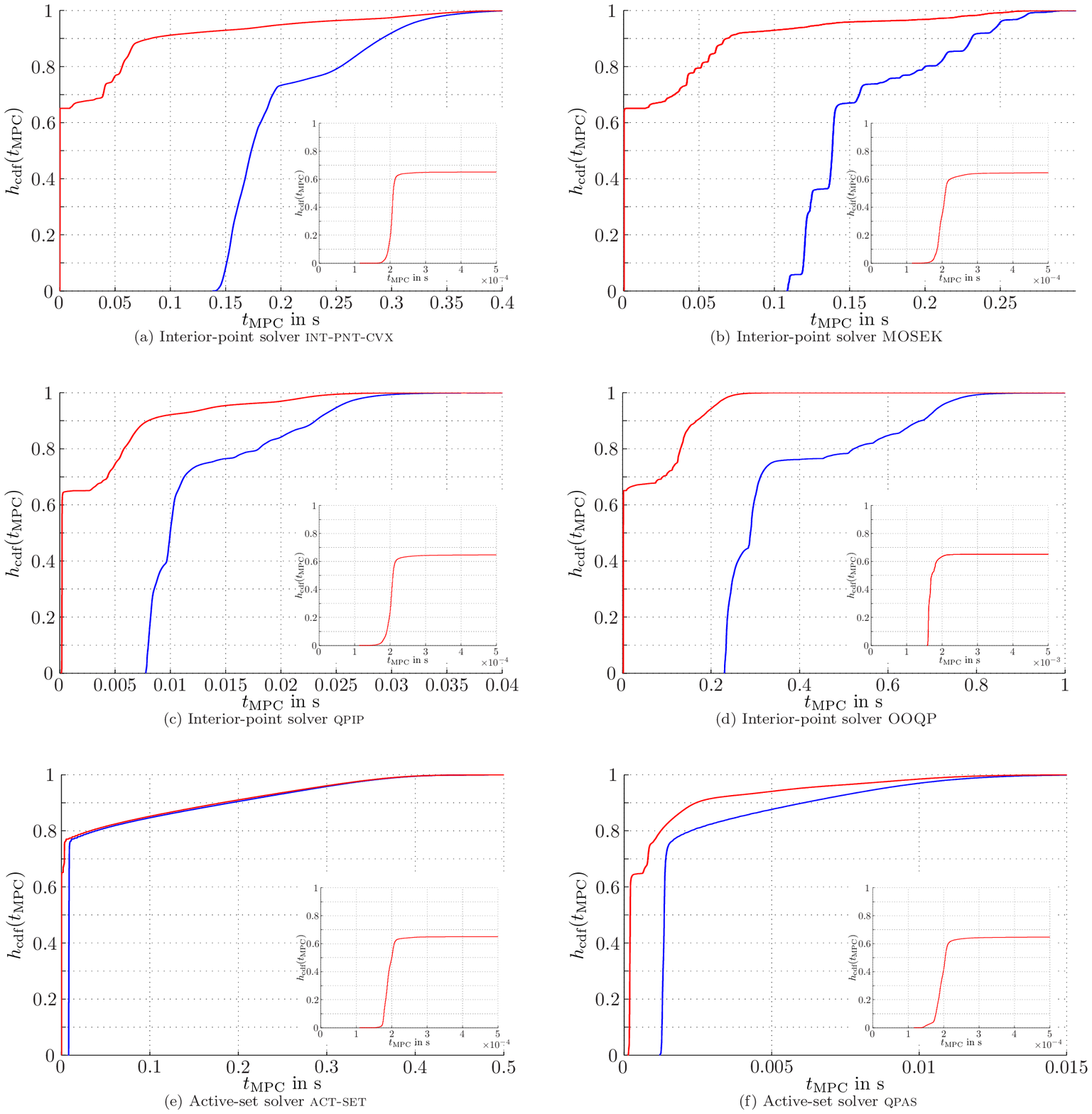}
  \caption{Cumulative distribution functions for the COMA40 example (full-MPC in blue, CR-MPC in red).}\label{fig:resCOMA}
\end{figure}

Constraint removal introduces some computational overhead, since inactive constraints have to be detected (which requires the comparison of two real numbers per constraint in every iteration of the QP solver) and the reduced problem~\eqref{eq:ReducedCondensedQP} must be set up. Our computational experiments indicate that this additional effort is outweighted by the savings in all but very small MPC problems. In fact, constraint removal does not result in a globally better cumulative distribution function of the computational times for the combination of the smallest example (ACC25) with the dual active-set solver qpas (see Fig.~\ref{fig:resACC}f). (Globally better cumulative distribution functions do result for all other example-solver combinations with constraint removal.) Note that even for the ACC25-qpas example-solver combination about 54\% of all QPs are solved faster with constraint removal than the fastest time obtained for all ACC25-qpas QPs without constraint removal, however.

The accelerations are also considerable in the following sense: For all 36 examples (resp.\ 34, 28 examples) more than 50\% (resp.\ 75\%, 90\%) of the QPs MPC with constraint removal required less time than the \textit{shortest time} achieved without constraint removal (see Tab.~\ref{tab:crFasterThanFastestFull}). 

\newpage 

\bibliographystyle{IEEEtran}

\begin{thebibliography}{10}
\providecommand{\url}[1]{#1}
\csname url@samestyle\endcsname
\providecommand{\newblock}{\relax}
\providecommand{\bibinfo}[2]{#2}
\providecommand{\BIBentrySTDinterwordspacing}{\spaceskip=0pt\relax}
\providecommand{\BIBentryALTinterwordstretchfactor}{4}
\providecommand{\BIBentryALTinterwordspacing}{\spaceskip=\fontdimen2\font plus
\BIBentryALTinterwordstretchfactor\fontdimen3\font minus
  \fontdimen4\font\relax}
\providecommand{\BIBforeignlanguage}[2]{{%
\expandafter\ifx\csname l@#1\endcsname\relax
\typeout{** WARNING: IEEEtran.bst: No hyphenation pattern has been}%
\typeout{** loaded for the language `#1'. Using the pattern for}%
\typeout{** the default language instead.}%
\else
\language=\csname l@#1\endcsname
\fi
#2}}
\providecommand{\BIBdecl}{\relax}
\BIBdecl

\bibitem{Jost2014}
\BIBentryALTinterwordspacing
M.~Jost, G.~Pannocchia, and M.~M\"onnigmann, ``{Online constraint removal:
  accelerating MPC with a Lyapunov function},'' \emph{Automatica}, vol.~57, pp.
  164--169, 2015. [Online]. Available:
  \url{dx.doi.org/10.1016/j.automatica.2015.04.014}
\BIBentrySTDinterwordspacing

\bibitem{Rossiter1998}
J.~A. Rossiter, G.~Kouvaritakis, and M.~J. Rice, ``A nummerically robust
  state-space approach to stable-predictive control stratigies,''
  \emph{Automatica}, vol.~34, pp. 65--73, 1998.

\bibitem{Boyd2009}
S.~Boyd and L.~Vandenberghe, \emph{Convex Optimization}.\hskip 1em plus 0.5em
  minus 0.4em\relax {Cambridge University Press}, 2009.

\bibitem{Tondel2003}
P.~T{\o}ndel, T.~A. Johansen, and A.~Bemporad, ``An algorithm for
  multi-parametric quadratic programming and explicit {MPC} solutions,''
  \emph{Automatica}, vol.~39, pp. 489 -- 497, 2003.

\bibitem{mpttoolbox}
M.~Kvasnica, P.~Grieder, and M.~Baoti\'{c}, ``{Multi-Parametric Toolbox
  (MPT)},'' 2004, {http://control.ee.ethz.ch/}.

\bibitem{Naus2010}
G.~J.~L. Naus, J.~Ploeg, M.~J.~G. Van~de Molengraft, W.~P. M.~H. Heemels, and
  M.~Steinbuch, ``Design and implementation of parameterized adaptive cruise
  control: An explicit model predictive control approach,'' \emph{Control Eng.
  Practice}, vol.~18, no.~8, pp. 882 -- 892, 2010.

\bibitem{Oliveri2011}
A.~Oliveri, G.~J.~L. Naus, M.~Storace, and W.~P. M.~H. Heemels,
  ``Low-complexity approximations of {PWA} functions: A case study on adaptive
  cruise control,'' in \emph{Proc. 20th Euro. Conf. Circuit Theory and Design},
  2011, pp. 669--672.

\bibitem{Lunze2004}
J.~Lunze, \emph{Regelungstechnik 2}.\hskip 1em plus 0.5em minus 0.4em\relax
  Springer Verlag, 2004.

\bibitem{Wang2010}
Y.~Wang and S.~Boyd, ``Fast model predictive control using online
  optimization,'' \emph{IEEE Transactions on Control Systems Technology},
  vol.~18, no.~2, pp. 267--278, 2010.

\bibitem{Wang2008}
\BIBentryALTinterwordspacing
------, ``fast\_mpc: code for fast model predictive control,'' 2008. [Online].
  Available: \url{http://stanford.edu/$\sim$ boyd/fast\_mpc/}
\BIBentrySTDinterwordspacing

\bibitem{OptToolBox}
\BIBentryALTinterwordspacing
MathWorks, \emph{MathWorks Matlab 2013a Documentation Center: Quadratic
  Programming Algorithms}. [Online]. Available:
  \url{http://www.mathworks.de/de/help/optim/ug/quadratic-programming-algorithms.html}
\BIBentrySTDinterwordspacing

\bibitem{qpc}
\BIBentryALTinterwordspacing
D.~A. Wills, \emph{QPC - Quadratic Programming in C}, School of Electrical
  Engineering and Computer Science, University of Newcastle, 2009-08-11.
  [Online]. Available: \url{http://sigpromu.org/quadprog/index.html}
\BIBentrySTDinterwordspacing

\bibitem{OOQP}
E.~M. Gertz and S.~J. Wright, \emph{OOQP User Guide}, october 2001, updated may
  2004~ed., Argonne National Laboratory, 9700 South Cass Avenue, Argonne, IL
  60439, 2004.

\bibitem{Gertz2003}
------, ``Object-oriented software for quadratic programming,'' \emph{ACM
  Transactions on Mathematical Software}, vol.~29, pp. 58--81, 2003.

\bibitem{MOSEK}
\BIBentryALTinterwordspacing
\emph{MOSEK version 7 manuals}, MOSEK ApS. [Online]. Available:
  \url{http://www.mosek.com/resources/doc}
\BIBentrySTDinterwordspacing

\bibitem{Nocedal2000}
J.~Nocedal and S.~J. Wright, \emph{Numerical Optimization}, 2nd. Ed., Ed.\hskip
  1em plus 0.5em minus 0.4em\relax Springer Verlag, 2000.

\bibitem{Goldfarb1983}
D.~Goldfarb and A.~Idnani, ``A numerical stable dual method for solving
  strictly convex quadratic programs,'' \emph{Mathematical Programming},
  vol.~27, pp. 1--33, 1983.

\end{thebibliography}

\end{document}